\documentclass[runningheads]{llncs}
\usepackage[ruled,vlined,algochapter,linesnumbered, boxed]{algorithm2e}
\usepackage{amsmath,amssymb,mathabx}
\usepackage[T1]{fontenc}
\usepackage{todonotes}
\usepackage{booktabs}
\usepackage{cleveref}
\usepackage{pdflscape}
\usepackage{afterpage}
\usepackage{capt-of}
\usepackage{xspace}
\usepackage{mathtools}
\usepackage{etoolbox}
\usepackage{multirow}
\usepackage{url}

\usepackage{charter}
\usepackage[scaled=0.90]{FiraMono}

\DeclarePairedDelimiter\abs{\lvert}{\rvert}

\newcommand{\integerx}[1]{\hat{x}^{(#1)}}
\newcommand{\fracx}[1]{\bar{x}^{(#1)}}
\newcommand{\fracy}[1]{\bar{y}^{(#1)}}

\newcommand{\Z}{\mathbb{Z}}

\newcommand{\scip}{SCIP{}}
\newcommand{\soplex}{SoPlex{}}
\newcommand{\fandp}{fix-and-propagate\xspace}
\newcommand{\Fandp}{Fix-and-propagate\xspace}
\newcommand{\norm}[1]{\left\lVert#1\right\rVert}
\newcommand{\loss}[1]{\textit{\textcolor{red}{#1}}}
\newcommand{\win}[1]{\textbf{\textcolor{blue}{#1}}}

\usepackage{graphicx}
\begin{document}
\title{Scylla: a matrix-free fix-propagate-and-project heuristic for mixed-integer optimization}
\titlerunning{Scylla: a matrix-free heuristic for mixed-integer optimization}
%
\author{
Gioni Mexi$^*$\inst{1}\orcidID{0000-0003-0964-9802}\and
Mathieu Besançon\inst{1}\orcidID{0000-0002-6284-3033} \and
Suresh Bolusani\inst{1}\orcidID{0000-0002-5735-3443} \and
Antonia Chmiela\inst{1}\orcidID{0000-0002-4809-2958} \and
Alexander Hoen\inst{1}\orcidID{0000-0003-1065-1651} \and
Ambros Gleixner\inst{1,2}\orcidID{0000-0003-0391-5903}
}
\authorrunning{G. Mexi et al.}
%
\institute{Interactive Optimization and Learning, Zuse Institute Berlin, Germany \and HTW Berlin, Germany\\
\email{\{mexi, besancon, bolusani, chmiela, hoen, gleixner\}@zib.de}
}

\let\oldmaketitle\maketitle
\renewcommand{\maketitle}{\oldmaketitle\setcounter{footnote}{0}}

\maketitle              
\begin{abstract}
We introduce Scylla, a primal heuristic for mixed-integer optimization problems.
It exploits approximate solves of the Linear Programming relaxations through the matrix-free Primal-Dual Hybrid Gradient algorithm with specialized termination criteria,
and derives integer-feasible solutions via fix-and-propagate procedures and feasibility-pump-like updates to the objective function.
Computational experiments show that the method is particularly suited to instances with hard linear relaxations. 

\end{abstract}

\section{Introduction}\label{sec:intro}
\renewcommand*{\thefootnote}{\fnsymbol{footnote}}
\footnotetext[1]{Corresponding author.}
\renewcommand*{\thefootnote}{\arabic{footnote}}

Mixed-integer optimization (MIP) is an essential and challenging class of problems
found in numerous applications including transport, power systems, engineering, or manufacturing.
Primal heuristics are used as stand-alone methods to compute good-quality solutions or as subroutines in
branch-and-bound solvers, potentially exploiting the information extracted in the rest of the solution process~\cite{Berthold2014}.
Solutions derived from primal heuristics can be returned immediately when the corresponding gap is small,
they help to prune nodes of the branch-and-bound tree, and are also leveraged by other subroutines
of the solver, e.g., cut selection~\cite{SCIP6}.

When integrated within a branch-and-bound scheme, many primal heuristics assume that an optimal solution of the linear programming (LP)
relaxation is readily available and that the LP can be re-optimized quickly after few modifications.
These repeated calls to LP solvers can become a bottleneck in several scenarios.
First, while it is true that standard MIP solvers produce and rely on optimal LP relaxation solutions,
they are not yet available during presolving.
Second, for some applications, re-optimization via the simplex method is expensive~\cite{KochEtAl2011}.
For primal heuristics in a stand-alone setting, the computation of an initial
solution to the LP relaxation by standard methods that require matrix factorizations,
such as the simplex or interior point methods, may be expensive for large-scale applications.

To address this problem, we present \emph{Scylla}, a \emph{matrix-free} primal heuristic, in the sense that it uses only iterated 
matrix-vector multiplications to compute approximate LP solutions without 
factorizing the constraint matrix nor directly solving linear systems.
We supplement this with LP-free constraint propagation to round such approximate reference solutions to integer feasibility and objective function updates inspired by the well-known feasibility pump~\cite{fischetti2005feasibility,bertacco2007feasibility,achterberg2007improving}.

The paper is organized as follows.
In \Cref{sec:algorithm}, we present the algorithm design and highlight its key features.
In \Cref{sec:comp-results}, we report on our computational experiments.
We summarize our observations and discuss future research directions in \Cref{sec:conclusions}.

\section{The algorithm}\label{sec:algorithm}
\newcommand{\pstep}{\tau_x}
\newcommand{\dstep}{\tau_y}
\newcommand{\pdir}{d_x}
\newcommand{\ddir}{d_y}

The overall structure of Scylla is outlined in \Cref{algo:overall}
and is based upon several subroutines:
\begin{itemize}
\item a limited-precision call to the Primal-Dual Hybrid Gradient (PDHG) algorithm~\cite{zhu2008efficient,chambolle2011first} to compute an approximate solution of the LP relaxation,
\item a \fandp procedure to attempt rounding this fractional solution to a fully, i.e., also integer, feasible solution,
\item an objective function update similar to that of the feasibility pump~\cite{fischetti2005feasibility,bertacco2007feasibility,achterberg2007improving} in order to iteratively drive the PDHG solution towards integer feasibility.
\end{itemize}
The outer loop of the algorithm is relatable to a feasibility pump~\cite{fischetti2005feasibility,bertacco2007feasibility} where the rounding step
uses constraint propagation~\cite{fischetti2009feasibility} and the projection step is replaced with a limited number of
PDHG iterations until relaxed termination criteria are met.
The inner iterations counted by~$k$ are dominated by matrix-vector multiplications with the constraint matrix and its transpose; the outer iterations counted by~$K$ are dominated by a \fandp call and update the objective function at the end.
Note that before starting the actual algorithm, we presolve the instance to 
reduce its size using the \scip{} solver~\cite{scip8}.

In the following sections, we describe the principal components of the algorithm and their interaction.

\begin{algorithm}[t]
	\DontPrintSemicolon
	\SetKwInOut{Input}{Input}\SetKwInOut{Output}{Output}
	\SetKwInOut{Init}{Initialization}
	\Input{MIP instance $\min\{c^\top x \, : \, Ax \geq b,\ \ell \leq x \leq u,\
   x_j \in \Z,\ \forall\ j \in I\}$}
	\Output{MIP-feasible solution $x$ or $\varnothing$ if none found}
	initialize $k = 1$, $K = 1$, $\alpha = 0.9$\;
    initialize primal-dual step sizes $\pstep > 0, \dstep > 0$ and $\fracx{k} = 0$, $\fracy{k} = 0$\;
	initialize $\hat{c} = c$\;
	\While{\textnormal{limits not reached}}{
		\Repeat(){$\fracx{k},\fracy{k} \approx$ LP optimal \text{or} bounds stalling \text{or} $k > \mathtt{iterlimit}$\label{line:pdhgloop2}}{\label{line:pdhgloop1}
			{
				update step sizes $\pstep,\dstep$\tcc*{PDHG iterations}
				$\pdir = A^\top \fracy{k} - \hat{c}$ and
				$\fracx{k+1} = \texttt{proj}_{[\ell,u]}(\fracx{k} + \pstep \pdir)$\;
				$\ddir = b - A(2 \fracx{k+1} - \fracx{k})$ and
				$\fracy{k+1} = \texttt{proj}_{\mathbb{R}^n_{+}}(\fracy{k} + \dstep \ddir)$\;
				  $k = k+1$\;
			  }\vspace*{1ex}
		  }
		\lIf{$\fracx{k}$ MIP-feasible}{\Return{$\fracx{k}$}\label{line:success1}}
		$\integerx{K} = \texttt{fix-and-propagate}(\fracx{k})$  \tcc*{rounding attempt}
		\lIf{$\integerx{K}$ MIP-feasible}{\Return{$\integerx{K}$}\label{line:success2}}
		\lIf{cycling on $\integerx{K}$ detected}{
			$\integerx{K} = \texttt{perturb}(\integerx{K}$)}\label{line:perturb}
		$\hat{c}^\top x = \alpha^{K}\frac{\sqrt{\abs{I}}}{\norm{c}} c^\top x + (1 - \alpha^{K})\Delta(x,\integerx{K}) $ \tcc*{objective update}\label{line:objup}
		$K = K+1$ 
	}
	\Return{$\varnothing$}
   \caption{High-level structure of the Scylla heuristic\label{algo:overall}}
\end{algorithm}

\subsection{Approximate LP solutions by matrix-free linear optimization}

Heuristic and exact methods for mixed-integer problems often rely on
repeated calls to a linear optimization solver.
These calls represent a substantial amount of the runtime and require, for the most common LP solving
techniques, a matrix factorization.
As highlighted in \cite{applegate2021practical}, first-order methods applied to the saddle-point formulations of LPs open the door
to a trade-off between solution accuracy and runtime which was rarely reasoned about in the linear and mixed-integer setting.
In our algorithm, these approximate LP solutions are used as a starting point for a \fandp procedure.

To leverage this trade-off, we use the PDHG \cite{chambolle2011first} algorithm and its
recent PDLP implementation
(PDHG specialized for linear optimization) \cite{applegate2021practical}.
Some notable algorithmic improvements were introduced in \cite{applegate2021practical} compared to the
standard PDHG including adaptive step size strategies, restarts, and matrix preconditioning.
We use the PDLP implementation available in OR-Tools\footnote{The package is available at \url{https://github.com/google/or-tools}.}.
The modifications allow us to warm-start the solution process from a previously computed primal-dual pair
and to add custom stopping criteria detailed in \Cref{sec:stopping}.

For conciseness of the algorithm description, we only included the core PDHG iterations and did not detail the step size computation procedure or restarting scheme from PDLP.
For Lines~\ref{line:pdhgloop1} to \ref{line:pdhgloop2} of 
Algorithm~\ref{algo:overall}, Scylla calls PDLP to solve similar LPs where only 
the objective function vector changes.
We use the fact that the state of PDLP is fully described by its primal-dual pair and primal and dual residuals to
warm-start all calls after the first one, which may help to reduce the total number of iterations after the first solve.
Using PDLP also means that Scylla remains a matrix-free heuristic only requiring sparse matrix-vector multiplications.

\subsection{Custom stopping criteria}\label{sec:stopping}

The standard stopping criterion for PDLP is a maximum error $\varepsilon$ on the
primal and dual feasibilities.
This error can be relaxed compared to what would be expected in ``standard'' LP solving since the resulting fractional solution
is then guiding a rounding procedure. In particular, this error can also be 
adapted dynamically to allow a bigger tolerance at the beginning of the 
heuristic.
In addition to increasing the initial error tolerance (to
$\varepsilon^{(1)}=0.01$ in our experiments), we also use a progressive
refinement scheme in which the error is reduced geometrically by a factor $\beta
= 0.98$, resulting in the following.
\begin{equation*}
	\varepsilon^{(K+1)} = \max\{\beta\ \varepsilon^{(K)}, 10^{-8}\}
\end{equation*}
\noindent

\subsection{The \fandp procedure}
\label{subsec::fix_and_propagate}

\Fandp algorithms have been used successfully as part of branch-and-bound solvers \cite{BertholdHendel2014,gamrath2019structure}.
They are also an essential step in the most successful variants of feasibility pump algorithms~\cite{fischetti2009feasibility}, and we use them with the same motivation: to ``round'' the PDHG solution~$\fracx{k}$ to integer feasibility.

``Rounding'' is achieved by fixing an integer variable with fractional value
to an integer within its domain, followed by a propagating step. This is 
repeated until all such variables are fixed or infeasibility is detected by
some domain becoming empty.
In the latter case,
\fandp continues in order to produce an integer vector by ignoring any constraint that would lead to empty domains.
At the end of the propagation, all remaining unfixed variables are fixed to their values in the fractional reference solution or its projection on to their domain.
The procedure always produces an integer-feasible, but not necessarily LP-feasible, solution.

One major issue of feasibility pump-like algorithms and hence also Scylla is cycling: after some iterations, it may happen that $\integerx{K} = \integerx{K'}$ for $1 \leq K' < K$. 
Similar to \cite {fischetti2005feasibility}, this issue is handled by randomly
perturbing some of the variables in the current rounded vector $\integerx{K}$
(Line~\ref{line:perturb}).

\subsection{Objective function update}

Inspired by objective feasibility pump approaches 
\cite{fischetti2005feasibility,bertacco2007feasibility}, we update the 
objective function as follows. After computing the integer point $\integerx{K}$ 
and before the subsequent PDLP call,
the objective function is updated (Line~\ref{line:objup} of
Algorithm~\ref{algo:overall}) as a convex combination of the initial objective function
$c^\top x$ and an $\ell_1$-norm distance to the integer
vector
$\Delta(x, \integerx{K})$
with $\alpha \in [0,1)$ chosen as constant $0.9$ in our experiments.

This convex combination represents a trade-off between the initial objective (towards high-quality solutions) and
the distance to an integer vector (towards integer feasibility).
After each iteration~$K$, the weight shifts further towards integer feasibility in the attempt to produce a PDLP solution $\fracx{k}$ that \fandp can successfully convert to a MIP-feasible solution.

\section{Computational results}\label{sec:comp-results}
As a baseline for our implementation, we modified an open source implementation of 
the Feasibility Pump~\cite{mexi2023using} 
(FP)\footnote{\url{https://github.com/GioniMexi/FeasPumpCollection}} (that uses \soplex ~6.0.2~\cite{scip8} as the underlying LP solver) and extended it to
support approximate LP solves with PDLP.

Our test set consists of 
the benchmark library MIPLIB 2017~\cite{gleixner2021miplib} and the instances used
in \cite{fischetti2009feasibility}. 
We excluded all instances that were solved to optimality during presolving with SCIP,
leaving us with 283 instances.

All results were obtained on a cluster of Intel Xeon Gold 6338 @ 2.0 GHz
machines with 1024 GB of RAM and a time limit of $3600$ seconds.
To mitigate the impact of performance variability~\cite{Lodi_2013}, both FP and Scylla are
tested on the complete test set with ten different random seeds.

\Cref{tab:results} shows the summarized results for 2830 instances (283 instances
$\times$ 10 random seeds). On the whole testset, FP
outperforms Scylla w.r.t.~both solution quality and running time. This is an 
expected result since the LPs in the FP iterations are for most instances 
easily solved to optimality, and warm-starting of simplex LP solvers is very effective.
Scylla therefore loses the benefit of exact LP solutions while still requiring 
similar costs per iteration. We see this when looking at the two different 
subsets of instances.

For the instances where the root LP solves in $[1,60]$
seconds, FP is able to find more solutions than Scylla and provides better 
quality solutions. Scylla is however faster on average, despite suffering under 
a greater model construction overhead which increases with the instance size.
Note that this overhead is tied to the current implementation and not to the 
heuristic itself. More specifically, the current OR-Tools interface through 
which PDLP is accessed requires re-building the model at every iteration of the 
feasibility pump, a restriction that penalized our approach, in particular for 
larger instances.

For instances in $[60, \text{lim}]$, Scylla significantly improves over FP in 
terms of both the number of solved instances and running time.
Notably, Scylla is capable of finding feasible solutions for 42 more instances 
than FP, and it is nearly four times faster.
These results suggest that for MIPs with difficult LP relaxations, Scylla is 
the preferred heuristic choice to compute first primal solutions.

 \begin{table}[t]
    \caption{Aggregated results comparing FP with Scylla. Row $[1,60]$ 
    ($[60,\text{lim}]$) only considers instances on which solving the 
    	LP takes between 1 and 60 seconds (more than 60 seconds). Column `\#' 
    	(`only') refers to the number of instances for which a solution was 
    	found (exclusively by only one heuristic), `gap' refers to the relative 
    	gap to the best known 
    	solution, `wins' refers to the number of instances with a performance
    	improvement of 10\% over the other heuristic, and `overhead' refers to
    	the overhead time to build the model.}
    \label{tab:results}
    \footnotesize
    \sffamily
    \setlength{\tabcolsep}{4pt}
    \centering
    \begin{tabular}{ll cc cc ccc}
        \toprule
          & & \multicolumn{2}{c}{found sol} & \multicolumn{2}{c}{\#wins} & 
          \multicolumn{2}{c}{shift. geom. mean}& \\
 \cmidrule(lr){3-4} \cmidrule(lr){5-6} \cmidrule(lr){7-8} 
              testset &algorithm &   \# &  only &           gap &       time(s) 
              &           gap &       time(s) & overhead(s)\\
		\midrule
	all		  &FP & 2268 &   152 &  \loss{1154 }&  \loss{1557} &  \loss{27.46} 
	&  \loss{10.06} & 0.02\\
			  &Scylla & \win{2270} &   \win{154} &  485 &  797 &  40.92 &   
			  17.30 & 1.93\\
			  \midrule
	$[1,60]$		  &FP & \loss{660}&   \loss{62} &    \loss{341} &       295 
	&  \loss{26.75} &  43.16 & 0.01 \\
			  &Scylla &  635 &    37 &  130 &  \win{381} &  41.78 &   
			  \win{36.75} & 1.99\\
			  \midrule
			  $[60, \text{lim}]$ &FP &  163 &    40 &             98 &      
			  64 &  \loss{28.11} &  974.66& 0.01\\
&Scylla & \win{ 205 } &   \win{ 82 } &  \win{105} &  \win{180} & 35.37 &    
\win{234.80} & 6.09\\
        \bottomrule
    \end{tabular}
 \end{table}

\section{Conclusion}\label{sec:conclusions}

In this paper, we presented Scylla, a primal heuristic for mixed-integer linear optimization
that combines matrix-free computation of approximate LP solutions with a \fandp procedure and an objective function update scheme similar to the feasibility pump.
Computational experiments on a large set of instances from the MIPLIB 2017
collection show that the method is particularly suited to instances with hard linear relaxations. 
Our current implementation does not exploit the potential for parallelization that lies in the sparse matrix-vector multiplications of PDHG, which dominate the PDLP runtime.
A robust, efficient, and multithreaded sparse BLAS library will thus accelerate this core operation of PDLP and could be key to a broader success especially on large-scale instances.

One aspect to be explored in future research is the fact that PDHG is by design suited to solve more general convex differentiable optimization problems, and PDLP handles quadratic problems
where the objective Hessian is diagonal \cite{pdlportools}.
On the one hand, this makes Scylla a candidate heuristic for the broader class of MIQPs.
On the other hand, this can be exploited also in the linear setting by choosing the $\ell_2$-norm as the distance function $\Delta(\cdot, \integerx{K})$ and avoiding auxiliary variables needed to model the $\ell_1$-norm on general integer variables \cite{fischetti2005feasibility,bertacco2007feasibility}.

\section*{Acknowledgements}

Research reported in this paper was partially supported through the Research Campus Modal funded by the German Federal Ministry of Education and Research (fund numbers 05M14ZAM,05M20ZBM) and the Deutsche Forschungsgemeinschaft (DFG) through the DFG Cluster of Excellence MATH+.

\bibliographystyle{splncs04}
\bibliography{refs}

\end{document}